\documentclass{amsart}
\usepackage{graphicx}
\usepackage{amsthm,amsfonts,amssymb}
\usepackage{amsmath}
\usepackage{hyperref}
\usepackage[latin1]{inputenc}
\usepackage{latexsym}
\newtheorem{theor}{Theorem}[section]
\newtheorem{lemma}{Lemma}
\newtheorem{proposition}{Proposition}
\newtheorem{corollaire}{Corollary}

\theoremstyle{definition}

\newcommand \be{\begin{eqnarray*}}
\newcommand \ee{\end{eqnarray*}}
\newcommand \ben{\begin{eqnarray}}
\newcommand \een{\end{eqnarray}}

\newcommand{\bigO}[1]{\mathcal O\pa{#1}}

\def \build#1#2#3{\mathrel{\mathop{\kern 0pt#1}\limits_{#2}^{#3}}}

\def \esp#1#2{\mathbb E_{#1}\left[#2\right]}

\def \build#1#2#3{\mathrel{\mathop{\kern 0pt#1}\limits_{#2}^{#3}}}

\newcommand{\ac}[1]{\left\{#1\right\}}

\newcommand{\pa}[1]{\left(#1\right)}

\newcommand{\pr}[1]{\mathbb{P}\left(#1\right)}

  \newdimen\AAdi%
\newbox\AAbo%
\def\AAk#1#2{\setbox\AAbo=\hbox{#2}\AAdi=\wd\AAbo\kern#1\AAdi{}}%
\title{Asymptotic analysis of a Drop-Push model
for Percolation}

\begin{document}

\author{Elahe Zohoorian Azad}
\address{Iran\\
Damghan\\
Damghan university \\
School of mathematics}

\email{ zohorian@dubs.ac.ir}


\begin{abstract}
In this article, we study a type of a one dimensional percolation
model whose basic features include a sequential dropping of
particles on a substrate followed by their transport via a pushing
mechanism (see [S. N. Majumdar and D. S. Dean, Phys. Rev. Ltt. A
11, 89 (2002)]). Consider an empty one dimensional lattice with
$n$ empty sites and periodic boundary conditions (as a necklace
with $n$ rings). Imagine then the particles which drop
sequentially on this lattice, uniformly at random on one of the
$n$ sites. Letting a site can settles at most one particle, if a
particle drops on an empty site, it stick there and otherwise the
particle moves according to a symmetric random walk until it takes
place in the first empty site it meet. We study here, the
asymptotic behavior of the arrangement of empty sites and of the
total displacement of all particles as well as the partial
displacement of some particles.
\end{abstract}
  \maketitle

%%%%%%%%%%%%%%%%%%%%%%%%%%%%%%%%%%%%%%%%%%%%%%%%%%%%%%%%%%%%%%%%%
\section{Introduction}
%%%%%%%%%%%%%%%%%%%%%%%%%%%%%%%%%%%%%%%%%%%%%%%%%%%%%%%%%

Fundamental in the domain of percolation is the manipulation of
dynamic sets: sets that can grow, shrink or otherwise change over
time. Some algorithms, like for example the Kruskal or Prim
algorithms, for the research of the \emph{minimal covering tree}
of a graph, involve the grouping of some distinct elements into a
collection of disjoint sets, and implementing two operations,
U{\small NION}, that unites two sets, F{\small IND} that finds
which set a given element belongs to, see \cite{CLR}.\\

In a basic model, clusters with different masses change, over
time, through space and when two clusters are sufficiently close
they merge into a single cluster, with a probability quantified,
in some sense, by a {rate kernel} $R$ depending on the masses, the
positions and the velocities of the two clusters. However, such a
model, including the spatial distribution of clusters and their
velocity, is still too complicated for analysis. A first
approximation was suggested independently by Marcus \cite{marcus}
and Lushnikov \cite{lush1,lush2}, considering kernels depending
only on the
masses of the clusters.\\

A {Marcus--Lushnikov process} \cite{Al} is a continuous-time
Markov process whose state space is the set of partitions of $n$
or, equivalently, the set of measures $\mu= \sum_k\
\frac{n(k,t)}n\ \delta_k$ with $\sum_kkn(k,t)=n$, on the set of
positive integers $\mathbb{N}$. The $k$'s stand for the sizes of
clusters and $n(k,t)$ is the number of clusters with size $k$ at
time $t$. The size--$k$ clusters provide a fraction
$\frac{k\,n(k,t)}n$ of the total size $n$. A Marcus--Lushnikov
process evolves by instantaneous jumps according to the rule
`\emph{{each pair of clusters} $(c_1,c_2)$ {merges
at rate} $R(c_1,c_2)/n$}', which $R$ is the {rate kernel} of the process.\\

%In other words, the system of the clusters jumps from the state
%$\mu$ to the state $\mu+\frac 1n (\delta_{x+y}-\delta_x-\delta_y)$
%at rate $R(x,y)/n$. More precisely, if at time $t$ the state of
%system is $(x_i)_{i\geq 1}$, the joint distribution of $t+T$, the
%time of the next merge, and of $(I,J)$, the indexes
%of the clusters which merge, can be written as follow:\\

%The \textit{additive Marcus--Lushnikov process} (with kernel
%$K(x,y)=x+y$) is embedded in Yao's spanning tree model and the
%classic parking model \cite{CMr}.

%Suppose that a family $(T_{i,j})_{1\leq i<j}$ of independent
%random variables, in $\mathbb R_{+}$, with following distribution,
%has given
%$$
%\mathbb P(T_{i,j}>t) = \exp(-K(x_i,x_j)t/n).$$ Then $T$ and
%$(I,J)$ are defined by the following relation:
%$$
%\inf_{1\leq i<j}T_{i,j} = T_{I,J} = T.
%$$
%We know then that $T_{I,J}$ and $(I,J)$ are independents, that
%$T_{I,J}$ follows an exponential law with parameter
%$\sum_{i,j}K(x_i,x_j)/n$, and that $$\mathbb
%P((I,J)=(i,j))=\frac{K(x_i,x_j)}{\sum_{k,\ell}K(x_{k},x_{\ell})}.$$

The sizes of trees in the forest of the spanning-tree model of
Yao, or the sizes of blocks of cars in the classic parking model
 \cite{CHN}, form an {additive} Marcus-Lushnikov
process in which the rate kernel is $R(c_1,c_2)=c_1+c_2$. The
sizes of the connected components of the random graph of
Erd\"os-Renyi \cite{ErdRen}, form a {multiplicative}
Marcus-Lushnikov process in which the rate kernel is
$R(c_1,c_2)=c_1 c_2$. The average costs of the Union-Find
algorithms, in the model of Erd\"os-Renyi, were studied by Knuth
\& Sch{\"o}nhage \cite{KS}, and Stepanov \cite{Stepanov}. In two
cases, the clusters are the connected components of a graph, and
the merging of two clusters is caused by the addition of an edge
between elements of these clusters. We can, in this article,
suppose that the initial state consists of $n$ cluster of size
$1$, which corresponds to a graph completely disconnected with $n$
vertices but any edge. There is $n-1$ merging between the initial
state monodisperse, $\delta_1$, and the final state, $\frac 1
n\delta_n$, of the Marcus-Lushnikov process. As it soon will be
seen, the model on which we work comprises the additive case: the
evolution of the sizes of clusters is described here by an
additive Marcus-Lushnikov
process.\\

In first analysis, we can distinguish three different regimes in
the evolution of the additive Marcus-Lushnikov processes,
according to the size $A^n_{k,1}$ of the largest cluster after the
$k$-th jump, with the interpretations concerning the fragmentation
of trees \cite{ADD,PAV} or the analysis of hashing algorithms
\cite{CHL}: the \emph{sparse regime} for the case in which if
$\sqrt n =o(n-k)$, $A^n_{k,1}/n$ tends to $0$ in probability; the
\emph{transition regime}, when $n-k=\Theta(\sqrt n)$, several
clusters of size $\Theta(n)$ coexist, and, once renormalized,
clusters' sizes converge to the
 widths of excursions of a stochastic processes related
to Brownian motion; and finally the \emph{almost full regime} for
the case if $n-k=o(\sqrt n)$, $A^n_{k,1}/n$ tends to $1$ in
probability, and a unique giant cluster of size $n-o(n)$
 coexists with smallest clusters with total size $o(n)$.

%%%%%%%%%%%%%%%%%%%%%%%%%%%%%%%%%%%%%%%%%%%%%%%%%%%
 \section{Main theorems }
\label{section:mergingcosts}
%%%%%%%%%%%%%%%%%%%%%%%%%%%%%%%%%%%%%%%%%%%%%%%%%%%

Considering a Marcus-Lushnikov process, at the $k$-th jump, two
clusters with respective sizes $(S_{k,n},s_{k,n})$, $S_{k,n}\ge
s_{k,n}$ are merged, at a cost that may depend on the sizes
$(S_{k,n}, s_{k,n})$. For instance, in some implementations, a
label is maintained for each element, signaling the set it belongs
to, and when merging two sets, one has to change the labels of the
elements of one of the 2 sets. Yao \cite{YAO}, Knuth \&
Sch{\"o}nhage \cite{KS}, studied two algorithms
\textit{Quick-Find} and \textit{Quick-Find-Weighted}. {Quick-Find}
updates the labels of one of the two sets, selected arbitrarily,
leading to a cumulated cost $C^{QF}_{n,m}=\sum_{k=1}^{m}B_{k,n},$
in which $B_{k,n}=S_{k,n}$ with probability $1/2$ and
$B_{k,n}=s_{k,n}$ with probability $1/2$. {Quick-Find-Weighted}
updates the smallest set at a cost of $c_{k,n}=s_{k,n}$, leading
to a cumulated cost $C^{QFW}_{n,m}=\sum_{k=1}^{m}s_{k,n}$. In
other contexts where coalescence of two sets occurs, costs of
interest are $L_{k,n}$, the size of one of the two sets chosen
randomly with a probability that is proportional to its size, i.e.
$L_{k,n}=S_{k,n}$ with probability $S_{k,n}/(S_{k,n}+s_{k,n})$ and
$L_{k,n}=s_{k,n}$ with probability $s_{k,n}/(S_{k,n}+s_{k,n})$.\\

In our model, where each particle while falling on an occupied
site moves according to a symmetric random walk until it finds an
empty site, the merging cost of two clusters of particles, at the
dropping moment of the $k$-th particle, is the \emph{movements} of
this particle in the cluster on which it falls, until it finds an
empty site. We indicate this movements by $M_{k,n}$, which is also
the necessary time for $k$-th particle to find an empty site.
Obviously, $M_{k,n}$ depends on the size of the corresponding
cluster (see the definition of the cluster in Section
\ref{embeddings}). The partial cumulated cost is then
$$C_{n,m}=\sum_{k=1}^{m}M_{k,n},$$  which is
interpreted as the total movements of the $m$ first particles
$(1\leq m\leq n)$. By the two following theorems, we study the
concentration of the partial cost $C_{n,\lceil\alpha n\rceil}$,
and the limit law of the total cost $C_{n,n-1}$,
when $n$ tends to infinity.\\

\begin{theor}
\label{coutpart_biased} For each $\eta\in (0,1)$, and each
$\varepsilon$ positive,
$$
\lim_n\pr{\sup_{\alpha \in [0,1-\eta]} \left|
 \frac{C_{n,\lceil\alpha n\rceil}}{n} - \ \frac{\alpha^2(\alpha^2-3\alpha+3)}
 {6(1 -\alpha)^3}\right|\ge\varepsilon}
 =
 0~.
$$
\end{theor}
\vspace{.5cm}
\begin{theor}
\label{tot displac} We have,
$$\frac{C_{n,n-1}}{n^{5/2}}\build{\longrightarrow }{}
{loi}\frac{\sqrt 2}{6} \xi,$$ where $\xi$ is a random variable in
which the distribution is characterized by its moments:
$$\mathbb E(\xi^k)=\frac{k!\sqrt\pi}{2^{(7k-2)/2}\Gamma(\frac{5k-1}{2})}
{\bar a}_k,$$ with
$${\bar a}_k=2(5k-6)(5k-4){\bar a}_{k-1}+\sum_{j=1}^{k-1}
{{\bar a}_j {\bar a}_{k-j}} \hspace{0,3cm} k\geq 2;\ {\bar a}_1=
\sqrt 2~.$$
\end{theor}

\vspace{0.5cm} The suite of this article is organized in the
following way: in Section \ref{embeddings}, we explain the
embedding of the additive Marcus-Lushnikov process in our model,
and we calculate the probability of the merging of two arbitrary
clusters. In Section \ref{moments}, using only the properties of
the symmetric random walk, we can calculate the two first moments
of the partial cost, $M_{k,n}$. Theorem \ref{coutpart_biased} is
proved in Section \ref{smolu}, thanks to the convergence of the
additive Marcus-Lushnikov process to the certain solutions of the
Smoluchowski equation, derived by the analytical arguments in
\cite{norris}; we use, more precisely, Theorem 3.1 of \cite{CHN}.
In Section \ref{asymptotic} we show that the cumulated cost of our
model can be approximated by an additive functional on Cayley
trees, induced by the tolls $(n^2)_{n\geq 1}$ (Proposition
\ref{lemme3}), which makes it possible to apply the results of
\cite{EZoh}.

%%%%%%%%%%%%%%%%%%%%%%%%%%%%%%%%%%%%%%%%%%%%%%%%%%%%%%%%%%%%%%%%%
\section{Embedding of the additive Marcus--Lushnikov process}
\label{embeddings}
%%%%%%%%%%%%%%%%%%%%%%%%%%%%%%%%%%%%%%%%%%%%%%%%%%%%%%%%%%%%%%%%%%%
We start with a description of the additive Marcus--Lushnikov
process that helps to understand its relations to the parking
scheme, generally: at the $k$-th step pick a first cluster $P$
with the probability $\frac{|P|}{n}$ among the $n-k+1$ clusters,
and call it the `predator' (being a size--biased pick, one
obtains, generally, a cluster larger than the average cluster);
then pick the `prey' $p$ uniformly among the $n-k$ remaining
clusters, and suppose that $P$ eat $p$, producing a unique cluster
with size $|P|+|p|$. Remark that if, alternatively, both clusters
are size--biased picks (resp. if both are uniform picks), we
obtain the \textit{multiplicative Marcus--Lushnikov process}
(resp. the \textit{constant
kernel Marcus--Lushnikov process}, also called Kingman's process).\\

Consider a lattice with $n$ sites at a circle, on which a set
$\mathcal P=\{1,\dots,n-1\}$ of $n-1$ particles drop successively
and eventually stick. Letting a site can settles only one
particle, each particle $p$ drops on a random site $f(p)$. If the
first chosen site $f(p)$ is on an empty site, the particle stick
there. On the other hand, if the site $f(p)$ is occupied, the
particle executes a symmetric simple random walk and finally it
stick on the first empty site which it meets. The first chosen
sites $\pa{f(p)}_{p\in \mathcal P}$ are assumed independent and
uniform on the $n$ sites, numbered from $1$ to $n$.\\

In this model, the clusters are formed by the occupied sites, with
the following conventions:
\begin{itemize}
\item there are as many cluster as there are empty sites, \item a
cluster contains an empty site and the set of consecutive
 occupied sites before (going clockwise) this empty site,
\item the size of the cluster is the number of sites constituting
it, including the empty site, \item if an empty site follows
another empty site, it is considered as a size--1 cluster of its
own.
\end{itemize}

\begin{figure}[tbp]
\centerline{\includegraphics[width=4.6cm]{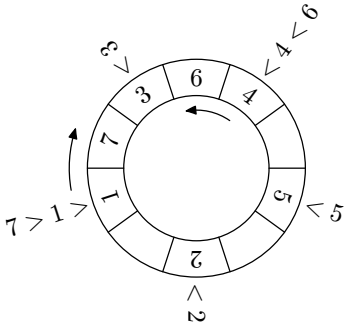} }\caption{ A
sample of tries $f(p)$ and the resulting 3 clusters.} Here
$n=10=6+2+2~.$
\end{figure}

\vspace{.4cm} The initial configuration, with $n$ empty sites, has
thus $n$ size--1 clusters (the \emph{monodisperse} configuration).
Each time that a particle sticks, two clusters merge, with
conservation of the mass, as the empty site that disappears and
the particle that replaces it both count for one mass--unit. The
final configuration, once the $n-1$ particles are sticken, is
constituted of a unique cluster with size $n$, and the unique
empty site, uniformly distributed on $\{1,2,\dots,n\}$. It turns
out that the sizes of clusters form an additive Marcus--Lushnikov
process, with kernel $K(x,y)=(x+y)/n$:
\begin{lemma}
\label{lemm:add M-L}Give that $k$ particles already sticken (that
$\ell=n-k$ sites are empty), and consider then two clusters with
sizes $x$ and $y$. The probability that these two clusters merge
at the next drop, $p_{n,k}(x,y)$, is
\begin{equation}
\label{pt} \frac{x+y}{n(n-k-1)}~.
\end{equation}
\end{lemma}

\begin{proof} Let $b_1,b_2,\dots,b_h;\ 1\leq h\leq k-1$ design
the non empty clusters just before the $k$-th particle drops (by
convention, $b_1$ design the cluster which contains the first
dropped particle). As the $k-1$ first particles choose uniformly
their sites, the order of these $h$ clusters, from $b_{1}$ on, is
a random uniform permutation. On the other side, it is not very
hard to see that the empty clusters (the clusters with size $1$)
are
merged uniformly on all their configurations:\\
Let $l_1,l_2,\dots,l_h$ being the number of clusters of size $1$
separating respectively the non empty clusters (i.e. the clusters
of size more than $1$) $b_1,b_2,\dots,b_h$. Thus dropping the
$k$-th particle, two empty clusters merge conditioning that these
two clusters be contiguous and the particle drops on the first
cluster (clockwise), in other words, conditioning that the $k$-th
particle drops on one of the
$$\sum_{i\in\{1,\dots,h\}}(l_i - 1)_{+}$$ empty sites surrounded
itself by two empty sites, on the right and on the left. The
conditional probability that two empty clusters merge, knowing the
position of the clusters, is thus
$$\frac 1n \sum_{i\in\{1,\dots,h\}}(l_i -
1)_{+}.$$ This probability does not depend on the sizes of the
nonempty clusters, but only on their number $h$, and the position
of the $r$ empty clusters among the $h$ nonempty clusters (note
that $k-1+h+r=n$). Let us pose now
 $r:=\sum_{i\in\{1,\dots,h\}}l_i$,
the full number of empty clusters, and let us note
 $\mathcal
L:=\{(l_1,l_2,\dots,l_h);\sum_{i\in\{1,\dots,h\}}l_i=r\}$, the set
of all configurations of $(l_1,l_2,\dots,l_h)$, often called
\textit{compositions} of $r$ with $h$ pieces. It is well-known
that
\[C^r_{h}=\text{Card } \mathcal L={{r+h-1}\choose{h-1}}.\]
Then the conditional probability that two empty clusters merge at
the $k$-th drop, knowing the number of empty and nonempty
clusters, is
\begin{eqnarray*}
p_{r,h,n} &=& \frac{\sum_{\mathcal L}\sum_{i\in\{1,\dots,h\}}(l_i
-1)_{+}}{nC^r_{h}}
\\
&=& \frac{1}{nC^r_{h}} \sum_{ \mathcal L}\left[r - h+\text{Card
}\{i;l_i=0\}\right]
\\
&=&\frac{1}{nC^r_{h}} \left[C^r_{h}(r-h)+\sum_{0\le\ell\le h-1}
\ell C^{r-h+\ell}_{h-\ell}{{h}\choose{\ell}}\right]
\\
&=&\frac{1}{nC^r_{h}} \left[C^r_{h}(r-h)+C^r_{h-1}h\right]
\\
&=& \frac{2{{r}\choose{2}}}{n(r+h-1)}~.
\end{eqnarray*}
Remark that $\sum_{0\le\ell\le h-1}\ell
C^{r-h+\ell}_{h-\ell}{{h}\choose{\ell}}$ can be interpreted as the
number of compositions of $r$ in $h$ pieces, such that a null
piece be marked, or underlined, whereas $C^r_ {h-1}h$ can be
interpreted as the number of compositions of $r$ in $h-1$ pieces,
such that one of the $h$ interstices between the pieces be marked,
or underlined: we obtain a bijective correspondence between two
sets inserting one zero additional into the site of the underlined
interstice, and underlining the zero so inserted. In addition, in
the additive Marcus-Lushnikov model, the probability that two
clusters of size 1 merge, at the stage $k$, if there is $r$ pieces
of size 1 and thus $h=n-k+1-r$ pieces of sizes higher than 1, is
also
\[{{r}\choose{2}}\ \frac{2}{n(n-k)}\] under the terms of Lemma
\ref{lemm:add M-L}.

Now let us consider the merging probability of two clusters of
respective sizes $x \geq 2$ and $y \geq 2$. Let us note $N_{x,y}$
the number of empty sites met while going clockwise, from the
cluster of size $x$ to the cluster of size $y$: $N_ {x,y}$ is
uniform on $\{1,2,\dots, n-k\}$. Obviously, if $N_{x,y}\notin \{1,
n-k\}$, the two clusters are not contiguous, and cannot merge
dropping of the $k$-th particle. For $N_{x,y}\in \{1,n-k\}$, let
us note $\delta$ the exit direction of cluster, $+$ or $-$
according to whether the particle leaves there in the clockwise
direction or in the opposite direction. Let us note $\tau$ the
size of the cluster in which the particle drops. We have then,
  \begin{eqnarray*}
  p_{n,k}(x,y)&=& \frac{1}{n-k}\sum_{\delta\in\{+,-\}}\sum_{\tau\in\{x,y\}}
 \sum_{N_{x,y}\in\{1,n-k\}}\mathbb{P}(\delta|\tau)\frac{\tau}n~.
 \end{eqnarray*}
 And as (see Section \ref{moments})
\begin{eqnarray*}
  \mathbb{P}(\delta|\tau) &=&
  \frac{\tau+1}{2\tau}\hspace{0,5cm}\text{if }\hspace{0,2cm}\delta=+,
 \\
 &=& \frac{\tau-1}{2\tau}\hspace{0,5cm}\text{if }\hspace{0,2cm}\delta=-,
 \end{eqnarray*}
we obtain well $p_{n,k}(x,y)=\frac{x+y}{n(n-k)}$. The merging
probability of a cluster of size $x\geq 2$ with one of the $r$
clusters of size 1 in this model, namely \[\frac{x+1}n\
\frac{C^{r-1}_{h}}{C^{r}_{h}},\] coincide also with the
probability in the additive Marcus-Lushnikov model, namely
$\frac{r(x+1)}{n(n-k)}$.
 \end{proof}
\vspace{0.4cm} From now, assume the Marcus--Lushnikov process to
be embedded in a drop particle scheme. In particular, we preserve
the interpretation of $L_{j,n}$ (resp. $R_{j,n}$) as size of the
cluster which is chosen by the $j$-th particle (resp. cluster
which merges with the preceding cluster when the $j$-th particle
sticks). We indicate by $p_{m,k}^{(j,n)}$ conditional probability
that $L_{j,n}$ (which we will interpret as the size of the $j$-th
predator before his meal) is equal to $k$ when the cluster created
by the dropping of $j$-th particle (the $j$-th predator after its
meal) is of size $m$. According to the asymptotic behavior of
$p_{m,k}^{(j,n)}$, when $m$ is large, we hope to reach a certain
intuition of the respective values of $L_{j,n}$ and $R_{j,n}$. It
proves, for combinative reasons, that $p_{m,k}^{(j,n)}$ do not
depend on $j$ or $n$: we have, for example,
\begin{eqnarray*}p_{m,k}^{(j,n)}=p_{m,k}^{(m-1,m)}&=&\mathbb P(L_{m-1,m}=k)\\
&=&\mathbb P(R_{m-1,m}=m-k)~.\end{eqnarray*} In what follows, we
shall remove thus the exponent of $p_{m,k}^{(j,n)}$.
\begin{lemma}
\label{lemm:comb-basic}
$$p_{m,k}={{m}\choose{k}} \frac{k^{k-1} (m-k)^{m-k-1}(2k-1)}{4 (m-1)m^{m-1}}
~.$$
\end{lemma}

\begin{proof} Let us calculate the probability $q_{m,k}^{(j,n)}$ that the
$j$-th merge utilizes a cluster of size $k$ and a cluster of size
$m-k$, knowing that the result of the merge is of size $m$: as it
is a probability concerning the evolution of the sizes of the
clusters, it is the same one for all the models where the
evolution of these sizes is described by an additive
Marcus-Lushnikov process. It is thus enough to calculate
$q_{m,k}^{(j,n)}$, as that was done in the parking model by
Chassaing and Marchand in \cite{CHN}. Here, we point out this
calculation for the convenience of the reader. Among the $n^j$
configurations for the $j$ first drops, there is
$${{j-1}\choose{m-2}}\ n\ m^{m-2} (n-m)^{j-m} (n-j-1)$$
configurations in which the $j$-th drop form a cluster of size
$m$: there is ${{j-1}\choose{m-2}}$ choice for the $m-2$ other
particles forming the cluster of size $m$, $n$ positions for this
cluster, and once the position and the particles are chosen, there
is $m^{m-2}$ ways to build this cluster of size $m$. The $j-m+1$
other particles can be sticken of $(n-m)^{j-m} (n-j -1)$ ways on
the $n-m-1$ sites which are reserved to them.

Among the configurations in which the $j$-th drop form a cluster
of size $m$, there is
\begin{eqnarray*}
{{j-1}\choose{k-1,m-k-1}}
&& n k^{k-2}(m-k)^{m-k-2}\\
 &&\times  (k+m-k)(n-m)^{j-m} (n-j-1)
 \end{eqnarray*}
configurations where the predator is of size $k$: there is ${{j-1}
\choose{k-1,m-k-1}}$ choice for the particles of the two clusters
intended to be merged, $n$ positions for this set of two adjacent
clusters, and once the particles of the two clusters and the
position are chosen, the $j-1$ first particles can stick in
$k^{k-2}(m-k)^{m-k-2}\ (n-m)^{j-m} (n-j-1)$ ways. This calculation
holds if we place the cluster of size $k$ initially, and there is
then $k$ choice for the site where the $j$-th particle drops. This
calculation holds also if we place the cluster of size $m-k$
initially, and there is then $m-k$ choice for the site where drops
$j$-th particle. It leads to
\begin{eqnarray*}
q_{m,k}^{(j,n)}&=& \frac{{{j-1}\choose{k-1,m-k-1}}k^{k-2}
 (m-k)^{m-k-2}\ m}{{{j-1}\choose{m-2}}m^{m-2}}
 \\
 &=& {{m}\choose{k}} \frac{k^{k-1} (m-k)^{m-k-1}}{(m-1)m^{m-2}}~.
\end{eqnarray*}

In our parking model, let $x$ denote the probability that there
exist two clusters of size $k$ and $m-k$ side by side, at the
dropping moment of the $j$-th particle. Then the probability that
the cluster of size $k$ is before (resp. after) the cluster of
size $m-k$ is $x/2$. If the cluster of size $k$ is on the left,
the predator is of size $k$ if the $j$-th particle falls on the
one of its $k-1$ occupied sites and exit from the right, with the
probability
\[\frac{x(k-1)}{4n},\] or if the $j$-th
particle falls on the single empty site of the cluster of size
$k$, with the probability \[\frac{x}{2n}.\] If the cluster of size
$k$ is on the right, the predator is of size $k$ if the $j$-th
particle falls on one of its $k-1$ occupied sites and exit from
the left, with the probability \[\frac{x(k-1)}{4n}.\] So the
cluster of size $k$ is the predator with the probability
\[\frac{xk}{2n},\] and the cluster of size $m-k$ is the predator
with the probability
\[\frac{x(m-k)}{2n}.\] We
deduced that
\[\frac{xm}{2n}=q_{m,k}^{(j,n)},\]
and that
\[p_{m,k}^{(j,n)}
= \frac{k}{m}\ q_{m,k}^{(j,n)}~.\] Finally
\begin{eqnarray*}
 p_{m,k}
 &=&
 {{m-1}\choose{k-1}} \frac{k^{k-1} (m-k)^{m-k-1}}{(m-1)m^{m-2}},
\end{eqnarray*}
as expected.\end{proof}

The Lemma \ref{lemm:comb-basic} and the Stirling formula entails
at once that
\begin{corollaire}
\label{rem_loideR}
\begin{equation} \label{convenloi}
\forall k \geq 1, \lim_{m \rightarrow \infty}
p_{m,m-k}=\frac{k^{k-1}e^{-k}}{k!}~.
\end{equation}
\end{corollaire}
Thus the limiting distribution of the size of the last prey is the
Borel distribution, in particular related to the explicit
solutions of Smoluchowski equations \cite{Al}, and to the function
of tree or Lambert function \cite{KNU2}. Thus, in law,
$R_{m-1,m}=\bigO{1}$. However, note that this distribution has an
infinite expectation, which is coherent with the fact that
$\esp{}{R_{m-1,m}}=\Theta\pa{\sqrt m}$. What we retain of these
calculations, is that provided $L_{k,n} + R_{k,n}$ be large,
$R_{k,n}$ or $s_{k,n}$ be negligible compared to $L_{k,n}$.

 \vspace{0.3cm}
%%%%%%%%%%%%%%%%%%%%%%%%%%%%%%%%%%%%%%%%%%%%%%%%%%
\section{Partial costs}
%\label{moments}
%%%%%%%%%%%%%%%%%%%%%%%%%%%%%%%%%%%%%%%%%%%%%%%%%%
In this section, firstly, we calculate the first and the second
moments of the {movements} of the $k$-th dropped particle,
$M_{k,n}$, which we will need for the proof of Theorem \ref{tot
displac}. Then, we deal with the first part of the demonstration
of Theorem \ref{coutpart_biased}, which is in fact a corollary of
Theorem 3.1 of \cite{CHN}, stated here as Theorem
\ref{coutpart_general}.

%%%%%%%%%%%%%%%%%%%%%%%%%%%%%%%%%%%%%%%%%%%%%%%%%%
\subsection{Moments}
\label{moments}
%%%%%%%%%%%%%%%%%%%%%%%%%%%%%%%%%%%%%%%%%%%%%%%%%%
As in the previous section, an arbitrary cluster of size $s$
consists of $s-1$ particles sticken successively and an empty site
in the $s$-th position. Consider a particle drops in this cluster,
on the one of these $s$ sites chosen randomly (in the uniform
way). This first choice, noted $\mathcal X_{0}$, is thus a uniform
random variable on $\{1,2,\dots,s\}$. Consequently,
\begin{eqnarray*}
\mathbb E(\mathcal X_{0}) = \sum_{i=1}^{s}i/s=\frac{s+1}{2}
\end{eqnarray*}
and
\begin{eqnarray*}
\mathbb E(\mathcal
X_{0}^2)=\sum_{i=1}^{s}{i^2}/s=\frac{(s+1)(2s+1)}{6}~.
\end{eqnarray*}
Consider now the variable
$$\mathcal X_{h}=\mathcal X_{0}+\sum_{i=1}^{h} Y_{i}\hspace{0,5cm};
\hspace{0,2cm}h=0,1,2,\dots,$$ representing the position of the
particle in the cluster of size $s$, after $h$ step. The $Y_{i}$
are the Bernoulli random variables of parameter $\frac{1}{2}$ with
value in $\{-1,1\}$. We indicate by $D_{s}$ the number of steps
inside the cluster of size $s$, before the particle sticks (that
one can also see as the time of receive to the edge of the
cluster). We have then
\begin{eqnarray*}
\mathbb{P}(\mathcal X_{D_{s}}=s) &=& \mathbb{P}(\mathcal
X_{D_{s}}=s|\mathcal X_{0}=s)\mathbb{P} (\mathcal X_{0}=s)\\
&+& \mathbb{P}(\mathcal X_{D_{s}}=s|\mathcal X_{0}\neq
s)\mathbb{P}(\mathcal X_{0}\neq s)
\\
&=& \frac{1}{s}+\frac{1}{2}\frac{s-1}{s}
\\
&=& \frac{s+1}{2s}~.
\end{eqnarray*}

Consequently, $\mathbb{P}(\mathcal X_{D_{s}}=0)=\frac{s-1}{2s}$.
Remark that $0$ indicates the last site before the cluster (in the
clockwise direction), site which is empty.\\
The processes $\mathcal M_{h}=\mathcal X_{h}^2-h$ and
\begin{eqnarray*}
 \mathcal M'_{h}=\mathcal X_{h}^4-2(3h-2)\mathcal
X_{h}^2+h(3h-1),
\end{eqnarray*}
are martingales. By stopping theorem,
$$\mathbb E(\mathcal X_{D_{s}}^2-D_{s})=\mathbb E(\mathcal X_{0}^2-0),$$
which gives
\begin{equation}
\label{esp T}\mathbb E(D_{s})=\frac{s^2-1}{6}~.
\end{equation}
Since the variables $\mathcal X_{D_{s}}$ and $D_{s}$ are not
independents, the calculate of the second moment of $D_s$ starting
from $\mathcal M'_{h}$ is not direct. To circumvent the
difficulty, we add the site $0$ to the cluster of size $s$, thus
obtaining a cluster of size $s+1$, such that we fall in a
symmetrical situation: we consider then an initial position
$\mathcal {\tilde X}_{0}$ uniform on $\{0,1,\dots,s\}$. We define
then $\tilde D_{s}$ and $\mathcal {\tilde X}_{\tilde D_{s}}$ in a
way similar to $D_{s}$ and $\mathcal {X}_{D_{s}}$, but these two
new variables are now independents. We have thus
\begin{eqnarray}
\label{T tildT} \mathbb E({\tilde D}_{s}^2) &=& \mathbb E({\tilde
D}_{s}^2|\mathcal{\tilde X}_{0}=0)\frac{s}{s+1}\nonumber\\
& +& \mathbb E({\tilde D}_{s}^2|\mathcal{\tilde X}_{0}\neq
0)\frac{s}{s+1}\nonumber\\
&=&\frac{s}{s+1}E({D}_{s}^2)~.
\end{eqnarray}
Moreover,
\begin{eqnarray*}
&&\mathbb E(\mathcal {\tilde X}_{0}^2) = \frac{s(2s+1)}{6},
\\
&&\mathbb E(\mathcal {\tilde
X}_{0}^4)=\frac{6s^5+15s^4+10s^3-s}{30(s+1)},
\\
&&\mathbb{P}(\mathcal {\tilde X}_{{\tilde D}_{s}}=s) = \frac12\ =
\ \mathbb{P}(\mathcal {\tilde X}_{{\tilde D}_{s}}=0)~.
\end{eqnarray*}
Stopping theorem for the martingale $$\tilde {\mathcal
M}_{h}=\mathcal {\tilde X}_{h}^4-2(3h-2)\mathcal {\tilde
X}_{h}^2+h(3h-1),$$ gives then
$$
\mathbb E[\mathcal {\tilde X}_{\tilde D_{s}}^4-2(3\tilde D_{s}-2)
\mathcal {\tilde X}_{\tilde D_{s}}^2+\tilde D_{s}(3\tilde
D_{s}-1)] = \mathbb E(\mathcal {\tilde X}_{0}^4+\mathcal {\tilde
X}_{0}^2),
$$
from which it is deduced that
\begin{equation}
\label{esp T^2} \mathbb E(D_{s}^2) = \frac{(s^2-1)(3s^2-7)}{45}~.
\end{equation}

%%%%%%%%%%%%%%%%%%%%%%%%%%%%%%%%%%%
\subsection{After $\lceil\alpha n\rceil$--th drop}
\label{smolu}
%%%%%%%%%%%%%%%%%%%%%%%%%%%%%%%%%%%%%
Here, Theorem \ref{coutpart_general}, gives the expression, in
terms of the solution $q(k,t)$ of Smoluchowski equation, of the
limit function $\varphi^\varsigma(\alpha)$ for the partial cost
\[
C^{\hat\varsigma}_{n,\lceil\alpha n\rceil} =
\sum_{k=1}^{\lceil\alpha
n\rceil}\hat\varsigma\pa{S_{k,n},s_{k,n},U_{k,n}},
\]
once $C^{\hat\varsigma}_{n,\lceil\alpha n\rceil}$ is normalized by
$n$. This theorem covers a wide class of costs, because the
general expression $\hat \varsigma\pa{S_{k,n},s_{k,n},U_{k,n}}$ of
the instantaneous cost of the $k$-th jump utilizes a randomization
parameter $U_{k,n}$, uniform on $[0,1]$. The asymptotic behavior
of the partial cost is expressed according to the conditional
instantaneous cost
\[
\varsigma(x,y) = \esp{}{
\left.\hat\varsigma\pa{S_{k,n},s_{k,n},U_{k,n}}\right\vert\pa{S_{k,n},s_{k,n}}
=(x,y) }~.
\]
Theorem \ref{coutpart_general} requires a hypothesis little
restrictive of polynomial growth of the moment of order 2 of the
instantaneous conditional cost,
\begin{equation}
\label{crpolyn} \forall x, y\in\mathbb{N},\ \
h(x,y)=\int_{0}^1\hat\varsigma\pa{x,y,u}^2\,du \le Ax^ny^m,
\end{equation}
for $A$, $m$ and $n$ well selected \footnote{In (\ref{crpolyn}),
$\mathbb{N}$ denotes the set of strictly positive entire
numbers.}. The cost $\hat{\varsigma}$ is supposed nonnegative, and
$(U_{k,n})_{k \in \mathbb{N}, n\in \mathbb{N}}$ denote a sequence
of independent random variables uniformly distributed on $[0,1]$.
We note $\varphi^\varsigma$ the increasing function of $[0,1)$ in
${\mathbb{R}}^+$ defined by
\[
\varphi^{\varsigma}(\alpha) = \int_0^{ \log \left(
     \frac{1}{1 -\alpha} \right)} \sum_{k \in \mathbb{N}} \sum_{l \in
  \mathbb{N}} \varsigma(k,l)\ (\frac{k+l}{2})\ q(k,t)q(l,t)\ dt,
\]
and
\[
q(k,t) = \frac{\left[k(1-e^{-t})\right]^{k-1} e^{-t}}{k!}\
\exp(-k(1-e^{-t}))~.
\]
We have thus
\begin{theor}[\cite{CHN}]
\label{coutpart_general} For all $\eta>0$,
$$\sup_{\alpha \in [0,1-\eta]}
\left|
 \frac{C^{\hat\varsigma}_{n,\lceil\alpha n\rceil}}{n} -
  \varphi^\varsigma(\alpha)
\right| \ \build{\longrightarrow }{}{P}\ 0~.$$
\end{theor}

We pose now some notations concerning our model: $p(k)$ denotes
the particle concerned with the $k$-th jump, the one which
verifies
\[
\#\ac{p\left\vert1\le p\le n-1\mbox{ and }T_p\le T_{p(k)}\right.}
= k,
\]
(where $T_p$ indicates the moment when the particle $p$ sticks).
Let us note $f_0(p(k))$, or $f_0(k)$ to be short, the first try of
$p(k)$, and note $f_j(k), j\geq 1,$ the $j$-th try of $p(k)$ at
the time of its search for an empty site. Let us note $\mathcal
H_k$ la $\sigma$-algebra generated by the trajectories
$\pa{f_j(\iota)}_{j\geq 0,\,1\le \iota\le k-1}$, the particles
$p(\iota);1\le \iota\le k-1$, and by $f_0(k)$. Conditioning by
$\mathcal H_k$,  we obtain, according to the relations (\ref{esp
T}) and (\ref{esp T^2}),
\begin{equation}
\label{rem_depl_cond1} \esp{}{M_{k,n}| \mathcal H_k} =
\frac{L_{k,n}^2-1}{6},
\end{equation}
and
\begin{equation}
\label{rem_depl_cond2} \esp{}{M_{k,n}^2| \mathcal H_k} =
\frac{(L_{k,n}^2-1)(3L_{k,n}^2-7)}{45}~.
\end{equation}

\vspace{0.3cm} The demonstration of Lemma \ref{lemm:comb-basic}
revealed that $L_{k,n}$ can be written
\[
L_{k,n} = S_{k,n}\mathbf{1}_{V_{k,n}\le
\frac{S_{k,n}}{s_{k,n}+S_{k,n}}} +s_{k,n}\mathbf{1}_{V_{k,n}>
\frac{S_{k,n}}{s_{k,n}+S_{k,n}}},
\]
where $V_{k,n}$ indicates randomly a number in $[0,1]$. To apply
Theorem \ref{coutpart_general}, we must write $M_{k,n}$ in the
form
\[
M_{k,n} = \hat\varsigma\pa{S_{k,n},s_{k,n},U_{k,n}}~.
\]

For that we must draw randomly $L_{k,n}$ in the set
$\ac{s_{k,n},S_{k,n}}$, using $V_{k,n}$, as explained above, then
we must randomly draw the first test of the $k$-th particle among
the $L_{k,n}$ sites of the clusters in which it falls, for example
in the form $\lceil L_{k,n}W_{k,n}\rceil$, where $W_{k,n}$
indicates randomly another number in $[0,1]$, independent of
$V_{k,n}$. Finally, it is necessary to simulate the random walk of
the $k$-th particle, for example using a sequence $\pa{Y_{k,n,
\ell}}_{\ell\ge 1}$ of independent random variables $\ac{\pm 1}$
symmetrical. It can be done, in a traditional way, by using the
coefficients of dyadic expansion

\[
U_{k,n} = \sum_{\ell\ge 1}\ \frac{d_{k,n, \ell}}{2^\ell}
\]
to reconstitute $\pa{V_{k,n}, W_{k,n},\pa{Y_{k,n, \ell}}_{\ell\ge
1}}$, like below
\begin{eqnarray*}
V_{k,n} &=& \sum_{\ell\ge 1}\ \frac{d_{k,n, 2\ell-1}}{2^\ell},
\\
W_{k,n} &=& \sum_{\ell\ge 1}\ \frac{d_{k,n, 4\ell-2}}{2^\ell},
\\
Y_{k,n,\ell} &=& 2d_{k,n, 4\ell}-1~.
\end{eqnarray*}
We have thus
\begin{eqnarray*}
\varsigma(x,y)&=&\esp{}{\left.M_{k,n}
\right\vert\pa{S_{k,n},s_{k,n}}=(x,y)}
\\
&=&\esp{}{\left.\esp{}{\left.M_{k,n} \right\vert {\mathcal H_k}}
\right\vert\pa{S_{k,n},s_{k,n}}=(x,y)}
\\
&=&\frac 16\  \esp{}{\left.\pa{L_{k,n}^2-1}
\right\vert\pa{S_{k,n},s_{k,n}}=(x,y)}
\\
&=& \frac 16\ \pa{\esp{}{ x^2\mathbf{1}_{V_{k,n}\le
\frac{x}{x+y}}+y^2\mathbf{1}_{V_{k,n}> \frac{x}{x+y}}} - 1}
\\
&=& \frac 16\ \pa{\frac{x^3+y^3}{x+y} - 1}~.
\end{eqnarray*}
In the same way,
\begin{eqnarray*}
h(x,y)&=&\esp{}{\left.M_{k,n}^2\right\vert\pa{S_{k,n},s_{k,n}}=(x,y)}
\\
&=&\frac{\esp{}{\left.(L_{k,n}^2-1)(3L_{k,n}^2-7)
\right\vert\pa{S_{k,n},s_{k,n}}=(x,y)}}{45}
\\
&=& \frac{(x^2-1)(3x^2-7)x+(y^2-1)(3y^2-7)y}{45(x+y)},
\end{eqnarray*}
such that $h$ be a polynomial and satisfies thus the assumption
\eqref{crpolyn}. We can then apply Theorem \ref{coutpart_general}
to $\hat \varsigma_{k,n}$ and $\varsigma(x,y)$~. Thus, for all
$\eta>0$,
$$\sup_{\alpha \in [0,1-\eta]}
\left|
 \frac{C_{n,\lceil\alpha n\rceil}}{n} - \varphi^\varsigma(\alpha)
\right| \ \build{\longrightarrow }{}{P}\ 0,$$ with
\begin{eqnarray*}
 &\varphi^{\varsigma}(\alpha)& = \\
&&\frac{1}{12} \int_0^{ \log \left(\frac{1}{1 -\alpha} \right)}
\sum_{k \in \mathbb{N}} \sum_{l \in\mathbb{N}} \pa{k^3+l^3-k-l}
q(k,t)q(l,t) dt
\\
&& =\frac{1}{6}\ \int_0^{ \log \left(\frac{1}{1 -\alpha} \right)}
 \langle q(.,t),x^3-x\rangle\;\langle q(.,t),1\rangle\ dt~.
\end{eqnarray*}
To finish calculation, we use the fact that $q(.,t)$ is expressed
in term of the Borel distribution, or of the Lambert function $T$
(cf. \cite{Jan2}). More precisely, if $B_{a}$ designs a Borel
random variable of parameter $a$, $0<a\le 1$, we have, for $k\ge
1$ : \[
 \pr{B_{a}=k} = \frac{(ka)^{k-1}}{k!}\ e^{-ka},
\]
and, for $a<1$,
\begin{eqnarray*}
&&\esp{}{B_{a}} = \frac{1}{1-a},\\
&&  \esp{}{B_{a}^2} = \frac{1}{(1-a)^3}, \\
&&\esp{}{B_{a}^3} = \frac{2a+1}{(1-a)^5}~.
\end{eqnarray*}
However it is noticed that, for $a=1-e^{-t}$,
\begin{eqnarray*}
q(k,t) &=& (1-a)\pr{B_{a}=k},
\\
\langle q(.,t),x^k\rangle & = & e^{-t}\esp{}{B_{a}^k},
\\
\langle q(.,t),x^3-x\rangle &=& 3e^{4t}-2e^{3t}-1,
\end{eqnarray*}
from where the calculation of $\varphi^{\varsigma}(\alpha)$
carried out higher, which leads to
\[
\varphi^{\varsigma}(\alpha) =
\frac{\alpha^2(\alpha^2-3\alpha+3)}{6(1 -\alpha)^3}~.
\]

\vspace{0.3cm}
%%%%%%%%%%%%%%%%%%%%%%%%%%%%%%%%%%%
\section{Cumulated cost}
\label{asymptotic}
%%%%%%%%%%%%%%%%%%%%%%%%%%%%%%%%%%%%%
To obtain Theorem \ref{tot displac}, it remains to estimate the
error made approximating $M_{k,n}$ by $(s_{k,n}+S_{k,n})^2$
(Proposition \ref{lemme3}), which makes it possible to reveal the
relation between $C_{n,n-1}$ and the additive functional with
penalties $(n^2)_{n\geq 0}$, studied in \cite{EZoh}. By keeping
this goal in memory, we demonstrate initially following lemmas. We
put
\[L_n = \sum_{k=1}^{n-1} L_{k,n}^2~.\]

\begin{lemma}
\label{lemme1} $\displaystyle \left\Vert
6C_{n,n-1}-L_{n}\right\Vert_2 = o\pa{n^{5/2}}$~.
\end{lemma}
\begin{proof} By developing $(6C_{n,n-1}-L_{n}+n-1)^2$,
we obtain:
\[
\left\Vert6C_{n,n-1}-L_{n}-n+1\right\Vert_2^2 = \Xi_1 +\Xi_2,
\]
where
\begin{eqnarray*}
\Xi_1 & = & \sum_{k=1}^{n-1}\ \esp{}{\pa{6M_{k,n}-L_{k,n}^2+1}^2
},\nonumber
\end{eqnarray*}
and
\begin{eqnarray*}
\Xi_2 = 2\sum_{1\le i<j\le n-1} \esp{}{(6M_{i,n}-L_{i,n}^2+1
)(6M_{j,n}-L_{j,n}^2+1)}. \nonumber
\end{eqnarray*}
In consequence of the relation (\ref{rem_depl_cond1}), for $i<j$,
\begin{eqnarray*}
\esp{}{ \esp{}{
\left(6M_{i,n}-L_{i,n}^2+1\right)\left(6M_{j,n}-L_{j,n}^2+1\right)
\left| \ \mathcal H_j \right.}} & = & 0,
\end{eqnarray*}
thus $\Xi_2$ disappears. According to (\ref{rem_depl_cond2}), we
have also
\begin{eqnarray*}
\esp{}{ \left.\left( 6M_{k,n}-L_{k,n}^2+1
 \right)^2 \right| L_{k,n}}
 =  \frac{7}{5}L_{k,n}^4-6L_{k,n}^2+\frac{23}{5}~.
\end{eqnarray*}
Therefor,
\begin{eqnarray*}
\Xi_1 &\leq&
  6\sum_{k=1}^{n-1} \esp{}{L_{k,n}^4}\\
& \le& 6 n^5 \int_0^1 \esp{}{\left( \frac{L_{\lceil\alpha
n\rceil,n}}{n}
    \right)^4 } d \alpha~.
 \label{equ_conv_int}
\end{eqnarray*}
\noindent According to \cite{PITT}, for $0< \alpha <1$,
 $( B^n_{\lceil\alpha n\rceil,1}/n)_{n \in \mathbb{N}}$
 ($B^n_{k,1}$ denotes the size of the longest
cluster after the $k$-th drop), converges in probability to $0$,
thus
\[
\lim_n\esp{}{ \pa{\frac{L_{\lceil\alpha n\rceil,n}}n}^4} = 0,
\]
and Lebesgue's dominated convergence Theorem completes the proof.
\end{proof}

\begin{lemma}
\label{lemme2} $\displaystyle
\left\Vert\sum_{k=1}^{n-1}(L_{k,n}+R_{k,n})^2-L_{n}\right\Vert_1 =
\bigO{n^2 \log n}$~.
\end{lemma}

\begin{proof}
\begin{eqnarray*}
\left\Vert\sum_{k=1}^{n-1}(L_{k,n}+R_{k,n})^2-L_{n}\right\Vert_1&=&
 \sum_{k=1}^{n-1}\esp{}{R_{k,n}^2}\\
 &+&2\sum_{k=1}^{n-1}\esp{}{R_{k,n}L_{k,n}}~.
 \end{eqnarray*}
Thanks to the Lemmas $4.2$ and $4.8$ of \cite{CHN}, for all
$k\in\{1,\dots,n-1\}$ we have
\begin{equation}
\label{R2} \sum_{k=1}^{n-1} \esp{}{R_{k,n}^2} = \bigO{n^2 \log n},
\end{equation}
 and
\begin{equation}
\label{RL}  \esp{}{R_{k,n}|L_{k,n} }=\frac{n-L_{k,n}}{n-k}~.
\end{equation}
Therefor,
\begin{eqnarray*}
n^{-2}\sum_{k=1}^{n-1}\esp{}{R_{k,n}L_{k,n}} &=& \int_0^1
\frac{n}{n-\lceil\alpha n\rceil}\esp{}{ \frac{L_{\lceil\alpha
n\rceil,n}}{n}
     } d \alpha\\
     &-&\int_0^1 \esp{}{ \frac{L_{\lceil\alpha
n\rceil,n}^2}{n(n-\lceil\alpha n\rceil)}
    } d \alpha~.
    \end{eqnarray*}
Since for $0< \alpha <1$,
 $( B^n_{\lceil\alpha n\rceil,1}/n)_{n \in \mathbb{N}}$ converges
 in
probability to $0$ (\cite{PITT}), we have
\[
\lim_n\esp{}{ \frac{L_{\lceil\alpha n\rceil,n}}n} = 0,
\] { and } \[ \lim_n\esp{}{
\frac{{L^2}_{\lceil\alpha n\rceil,n}}{n(n-\lceil\alpha n\rceil)}}
= 0,
\]
and Lebesgue's dominated convergence Theorem
 completes the proof.
\end{proof}

$\ $

The two last lemmas involve the following proposition:

\begin{proposition}
\label{lemme3} \[\displaystyle \left\Vert
6C_{n,n-1}-\sum_{k=1}^{n-1}(L_{k,n}+R_{k,n})^2\right\Vert_1 =
o\pa{n^{5/2}}.\]
\end{proposition}

Now, $\sum_{k=1}^{n-1}(L_{k,n}+R_{k,n})^2$ is precisely the
additive functional on the Cayley trees induced by the penalties
$(n^2)_{n\geq 0}$ studied in \cite{EZoh}, to which the reader is
referred for more details. Though, we represent here Theorem 1.1
of \cite{EZoh} as the following Proposition:

\begin{proposition}
\label{Xn} Let $X_n$ be the additive functional defined on the
Cayley trees, induced by the toll $(n^2)_{n\geq 0}$. Then,
$$n^{-5/2}\ X_n \build{\longrightarrow }{}{\mathcal L}\ {\sqrt 2}\ \xi,$$
where $\xi$ is a random variable whose distribution is
characterized by its moments.:
$$\mathbb{E}(\xi^k)=\frac{k!\sqrt\pi}{2^{(7k-2)/2}\Gamma(\frac{5k-1}{2})}
{\bar a}_k,$$ where
$${\bar a}_k=2(5k-6)(5k-4){\bar a}_{k-1}+\sum_{j=1}^{k-1}
{{\bar a}_j {\bar a}_{k-j}} ;\hspace{0,1cm} k\geq 2\ ,\ {\bar
a}_1= \sqrt 2.$$
\end{proposition}

Finally Theorem \ref{tot displac} rises from Propositions
\ref{lemme3} and \ref{Xn}, thanks to the following theorem
\cite{Bil}:

\begin{theor}
\label{thecomb} Be $(X_n)_{n \in \mathbb{N}}$, $(Y_n)_{n \in
\mathbb{N}}$ and $X$ a random variable, such as for all $n$, $X_n$
and $Y_n$ are defined on the same probability space. If $(X_n)_{n
\in \mathbb{N}}$ converges in law to $X$ and if $(\Vert
X_n-Y_n\Vert )_{n \in \mathbb{N}}$ converges in probability to $0$
then $(Y_n)_{n\in\mathbb{N}}$ converges in law to $X$.
\end{theor}

\section*{Acknowledgments}
I wish to thank Philippe Chassaing for many fruitful discussions
that helped me to solve this problem.

\bibliographystyle{amsalpha}
\bibliography{TheseBib}

\providecommand{\bysame}{\leavevmode\hbox to3em{\hrulefill}\thinspace}
\providecommand{\MR}{\relax\ifhmode\unskip\space\fi MR }
\providecommand{\MRhref}[2]{%
  \href{http://www.ams.org/mathscinet-getitem?mr=#1}{#2}
}
\providecommand{\href}[2]{#2}
\begin{thebibliography}{CLR90}

\bibitem[Ald99]{Al}
D.J. Aldous, \emph{\it{Deterministic and stochastic models for coalescence
  (aggregation and coagulation): a review of the mean-field theory for
  probabilists}}, Bernoulli \textbf{5} (1999), 3--48.

\bibitem[AP98]{ADD}
D.J. Aldous and J.~Pitman, \emph{\it{The standard additive coalescent}}, Ann.
  Probab. \textbf{26} (1998), 1703--1726.

\bibitem[Bil95]{Bil}
P.~Billingsley, \emph{\it{Probability and measure}}, John Wiley \& Sons, New
  York, 1995.

\bibitem[CL02]{CHL}
P.~Chassaing and G.~Louchard, \emph{\it {Phase transition for parking blocks,
  Brownian excursion and coalescence}}, Random Structures Algorithms
  \textbf{21, no. 1} (2002), 76--119.

\bibitem[CLR90]{CLR}
T.H. Cormen, C.~E. Leiserson, and R.~L. Rivest, \emph{\it{Introduction to
  algorithms}}, McGraw-Hill, 1990.

\bibitem[CM04]{CHN}
P.~Chassaing and R.~Marchand, \emph{\it {Merging costs for the additive
  Marcus-Lushnikov process, and union-find algorithms}}, arxiv (2004),
  math.PR/0406094.

\bibitem[ER60]{ErdRen}
P.~Erd{\"o}s and A.~Renyi, \emph{\it {On the evolution of random graphs}},
  Publ. Math. Inst. Hungarian Acad. Sci. \textbf{Ser. A, 5A-2} (1960), 17--61.

\bibitem[Jan01]{Jan2}
S.~Janson, \emph{\it {Asymptotic distribution for the cost of linear probing
  hashing}}, Random Structures Algorithms \textbf{19(3-4)} (2001), 438--471.

\bibitem[Knu98]{KNU2}
D.E. Knuth, \emph{\it {Linear probing and graphs}}, Algoritmica \textbf{22,
  no.4} (1998), 561--568.

\bibitem[KS78]{KS}
D.E. Knuth and A.~Sch{\"o}nhage, \emph{\it {The expected linearity of a simple
  equivalence algorithm}}, Theoret. Comput. Sci. \textbf{6, no. 3} (1978),
  281--315.

\bibitem[Lus73]{lush1}
A.A. Lushnikov, \emph{\it {Evolution of coagulating systems}}, J. Colloid
  Interface Sci. \textbf{45} (1973), 549--556.

\bibitem[Lus78]{lush2}
\bysame, \emph{\it {Coagulation in finite systems}}, J. Colloid Interface Sci.
  \textbf{65} (1978), 276--285.

\bibitem[Mar68]{marcus}
A.H. Marcus, \emph{\it {Stochastic coalescence}}, Technometrics \textbf{10}
  (1968), 133--143.

\bibitem[Nor99]{norris}
J.R. Norris, \emph{\it {Smoluchowski's coagulation equation: uniqueness,
  nonuniqueness and a hydrodynamic limit for the stochastic coalescent}}, Ann.
  Appl. Probab. \textbf{9} (1999), 78--109.

\bibitem[Pav77]{PAV}
Yu.~L. Pavlov, \emph{\it {The asymptotic distribution of maximum tree size in a
  random forest}}, Th. Probab. Appl. \textbf{22} (1977), 509--5203.

\bibitem[Pit87]{PITT}
B.~Pittel, \emph{\it {Linear probing: the probable largest search time grows
  logarithmically with the number of records}}, J. Algorithms \textbf{8, no. 2}
  (1987), 236--249.

\bibitem[Ste70]{Stepanov}
J.~V.~E. Stepanov, \emph{\it {The probability of the connectedness of a random
  graph ${\mathcal G}\sb{m}\,(t)$}}, Teor. Verojatnost. i Primenen \textbf{15}
  (1970), 58--68.

\bibitem[Yao76]{YAO}
A.~C.~C. Yao, \emph{\it {On the average behavior of set merging algorithms}},
  Eighth Annual ACM Symposium on Theory of Computing (Hershey, Pa., 1976),
  Assoc. Comput. Mach., New York (1976), 192--195.

\bibitem[ZA]{EZoh}
E.~Zohoorian-Azad, \emph{\it {Limit Law of an additive functional on Cayley
  trees}}, to apear SIAM Journal on Discrete Mathematics.

\end{thebibliography}

\end{document}